\documentclass[preprint,12pt]{elsarticle}
\usepackage{amsmath,amsfonts,amssymb,dsfont,mathrsfs}
\usepackage{amssymb}
\usepackage{latexsym,stmaryrd}
\usepackage{pstricks}
\usepackage{textcomp}
\usepackage{graphicx}

\def\Ree{\mbox{$\mathds{R}$}}
\def\Iee{\mbox{$\mathds{I}$}}
\def\Cee{\mbox{$\mathds{C}$}}
\def\Nee{\mbox{$\mathds{N}$}}

\def\retard{\mbox{$\mathrm{e}^{-\vartheta s}$}}

\usepackage{hyperref}

\journal{Systems \& Control Letters}

\begin{document}

\begin{frontmatter}
\title{Approximation of distributed delays}
\author{Hao Lu\fnref{CSCship}}
\ead{hao.lu@insa-lyon.fr}
\author{Michael Di Loreto\corref{CorresAuthor}\fnref{CSCship}}
\ead{michael.di-loreto@insa-lyon.fr}
\author{Damien Eberard}
\ead{damien.eberard@insa-lyon.fr}
\author{Jean-Pierre Simon}
\ead{jean-pierre.simon@insa-lyon.fr}
\address{Laboratoire Amp\`ere, UMR CNRS 5005, INSA-Lyon, 20 Avenue Albert Einstein, 69621 Villeurbanne, France}
\cortext[CorresAuthor]{Corresponding author}
\fntext[CSCship]{The authors acknowledge the financial support of the CSC Scholarship from Xi'an Northwestern Polytechnical University (China), and the scholarship program with INSA-Lyon (France).}

\begin{abstract}
We address in this paper the approximation problem of distributed delays. Such elements are convolution operators with kernel having bounded support, and appear in the control of time-delay systems. From the rich literature on this topic, we propose a general methodology to achieve such an approximation. For this, we enclose the approximation problem in the graph topology, and work with the norm defined over the convolution Banach algebra. The class of rational approximates is described, and a constructive approximation is proposed. Analysis in time and frequency domains is provided. This methodology is illustrated on the stabilization control problem, for which simulations results show the effectiveness of the proposed methodology.
\end{abstract}

\end{frontmatter}

%====================================================================
\section{Introduction}\label{section1}
The interest for the use of distributed delays in the stabilization of time-delay systems appears in the pioneering work of~\citet{olbrot}. To generalize algebraic methods issued from linear systems in finite dimensional spaces to time-delay systems,~\citet{kamenetal} first introduce a general mathematical setting for the control, and in particular for the stabilization, of time-delay systems. This mathematical framework was formalized in~\cite{bretheloiseauFruit} by the introduction of the B\'ezout ring of pseudopolynomials, or independently by~\cite{glusingluerssen} with a behavioral approach. In these works, distributed delays appear in the stabilization problem. More generally, they are at the core in spectrum assignment by feedback~\cite{manitiusolbrot}. Finite spectrum assignment generalizes the principle of Smith's predictor for dead-time systems~\cite{smith} to more general systems with delays, which can be stable or unstable. Distributed delays appear also in the characterization of equivalence transformations~\cite{artstein}. Robustness issue and optimization~\cite{dymgeorgiousmithOptDelay}, robustness for input-delay systems~\cite{mondieniculesculoiseau} or finite time control~\cite{diloreto} are other topics where distributed delays play a central role.\medskip\\
A distributed delay is a linear input-output convolution operator of the form
\begin{equation}\label{eqintro1}
y(t)=(f\ast u)(t)=\int_{0}^{\vartheta}{f(\tau)u(t-\tau)\,\mathrm{d}\tau}
\end{equation} 
where $\vartheta$ is a strictly positive and finite real, and kernel $f(\cdot)$ is a continuous function with support $[0,\vartheta]$. Numerical implementation of distributed delay was early investigated, to propose effective algorithms for control. Such an implementation was the starting point of a widely research activity. A first proposition for approximation with finite dimensional systems was proposed in~\cite{kamenetalfinite}. Reduction and approximation of delay systems, involving lumped delays, were also investigated in~\cite{partington2004}. In the work of~\cite{vanassche}, the authors propose a numerical integral approximation to realize an operator like in~(\ref{eqintro1}). Such an approximation writes as a sum of lumped delayed distributions, and unfortunately introduces additional closed-loop poles, and also instability phenomena. See, \emph{e.g.}~\cite{zhongphd},~\cite{santosmondie}, and references therein. To overcome this problem, various solutions were proposed.
In~\cite{mirkinApproximation}, it was outlined that such an approximation has a poor accuracy and an high sensitivity for high frequencies. Hence, the author proposed to add a low-pass filter in the integral approximation. Such a solution was also proposed indepedently by~\cite{mondiemichiels}. Further implementation improvements were proposed in~\cite{zhong2005}, with rational approximation and extension of bilinear transformations. These last papers give positive answers to the open problem of general approximation of distributed delays, outlined in~\cite{richard}.\medskip\\
We propose in this paper a general methodology for numerical implementation of distributed delays. The numerical implementation of an operator exhibit two sides. The first one involves time discretization of the input-output behavior of the operator. Any distributed delay is a BIBO-stable operator. It is easy to verify that an equivalent discrete time system can be obtained by usual sampling tools, and that this system can be put into a sum of causal lumped distributions. Taking an appropriate sampling period, this equivalent discrete time system can always be made BIBO-stable. Hence, this first part presents no difficulty. We refer for instance to~\cite{zhong2004} and in references therein for more details on this part. The second part involves the approximation problem in continuous time of such an operator. This is this part we address in this paper. A continuous time approximation need to reproduce with high fidelity the internal dynamics of this operator, for large classes of input signals, but also to generate an arbitrarily close input-output behavior to the original one. For linear systems governed by a convolution, these notions are equivalent to impose simultaneous time and frequency approximations. For this, we explicitly separate two notions, namely input-output approximation and kernel approximation. We will explain why input-output approximation is not suitable for approximation of distributed delays, and we will then focus ourselves on kernel approximation. Kernel approximation can be realized in many different ways. Among most used classes of operators for approximation, we can cite polynomials, rational fractions, or exponentials. See, \emph{e.g.}~\cite{achieser},~\cite{cheney}, or~\cite{kammler}.\medskip\\ 
With the objective to substitute the distributed delay by a more tractable system, highlighting rational assumptions, we propose two classes to realize approximation, namely lumped systems and a subclass of distributed delays. With previous objectives, we enclose the approximation problem into the Wiener algebra of BIBO-stable systems, using the graph topology. This corresponds to the weakest topology where feedback is a robust property. Moreover, for stable systems, graph topology and norm topology being the same, we work on norm convergence over this algebra, which is a Banach algebra. This general framework was used first in~\cite{vidyasagaranderson} for approximation of distributed parameter systems by lumped systems. Roughly speaking, working over this algebra, an approximation comes down to realize approximation of the kernel over the Banach algebra $\mathscr{L}_1(\Ree_+)$. This idea also grew in~\cite{ohtaetal}, for the approximation of lumped delayed distributions which appear in optimal control. Here, we propose an extension of the classes of approximates, and we show that working in this general setting yields to approximate a distributed delay in both time and frequency domains, for large classes of input signals. We also propose to highlight the propositions made in~\cite{mirkinApproximation} or~\cite{mondiemichiels}, and to bring a mathematical foundation for this solution.\medskip\\
The paper is organized as follows. In Section~\ref{section2}, we define and fully characterize the main properties, in both time and frequency domains, of distributed delays. We characterize in particular a general decomposition of distributed delays on the so-called elementary distributed delays.
In Section~\ref{section3}, we explicit our approximation problem, and solve it. This section starts with general comments on approximation, where graph topology is briefly recalled, and where we show that any distributed delay can be approximated over this topology by lumped systems. Then, we move to another class of approximates, using a subclass of distributed delays, easy to implement with stability. We show a density property of this subclass, and analyze the convergence of the approximation which is proposed. Section~\ref{section4} relates the properties of such an approximation. We outline a proposal for a constructive approximation, and analyze time and frequency properties. A few simulation show the effectiveness of the method, on the stabilization control problem.
%====================================================================
\section{Convolution operators and distributed delays}\label{section2}
%====================================================================
\subsection{Convolution algebra}\label{section2-1}
An input-output causal convolution system is a dynamical system described by an equation of the form
\begin{equation}\label{eq2-1-1}
y(t)=(f\ast u)(t) = \int_{0}^{t}{f(\tau)u(t-\tau)\,\mathrm{d}\tau}
\end{equation}
where $y(\cdot)$, $u(\cdot)$ and $f(\cdot)$ are said to be the output, input and kernel of the map, respectively. Convolution systems are naturally defined over a commutative algebra since they are closed under addition, multiplication and scalar multiplication, operations that correspond to arbitrary series and parallel interconnections of such systems. For a normed algebra, it is quite interesting to obtain some closure properties for convergence. Hence, we prefer to work on a Banach algebra, that is an algebra for which any convergent sequence of elements in the algebra has a limit in the algebra, and for which the norm has the multiplicativity property. It will be assumed that all linear spaces and algebras are over the complex field. A general algebra of distributions including a wide class of convolution systems is given by the so-called Callier and Desoer algebra, denoted $\mathscr{A}$ \cite{callierdesoer}. We say that $f\in\mathscr{A}$ if 
\begin{equation}\label{eq2-1-2}
f(t) = \left\{
\begin{array}{ll}
f_{a}(t)+f_{pa}(t), & t\geq 0\\
0, & t<0
\end{array}\right.
\end{equation}
where the complex-valued function $f_{a}(\cdot)\in\mathscr{L}_{1}(\Ree_{+})$, that is $f_a$ is a complex valued function, locally integrable on $\Ree_+$, and such that $\int_{0}^{\infty}{|f_a(t)|\,\mathrm{d}t}<\infty$. The complex-valued distribution $f_{pa}$ stands for the purely atomic part and writes 
\begin{equation}\label{eq2-1-3}
f_{pa}(t) = \sum_{n=0}^{\infty}{f_{n}\delta(t-t_{n})},
\end{equation}
with $f_{n}\in\Cee$, $n=0,1, \ldots$, $0=t_0<t_1<t_2<\ldots$, $\delta(t-t_n)$ denotes the Dirac delta distribution centered in $t_n$, and $\sum_{n\geq 0}{|f_n|}<\infty$. As shown in \citet{desoervidyasagar}, it is well known that $\mathscr{A}$ is a commutative convolution Banach algebra with norm defined by
\begin{equation}\label{eq2-1-4}
\|f\|_\mathscr{A} = \|f_a\|_{\mathscr{L}_1}+\sum_{n=0}^{\infty}{|f_n|},
\end{equation}
and with unit element the Dirac delta distribution $\delta$. Denoting $\hat{f}$ the Laplace transform of $f$, $\hat{\mathscr{A}}$ denotes the set of Laplace transforms of elements in $\mathscr{A}$. The set $\hat{\mathscr{A}}$ is also a commutative Banach algebra with unit element under pointwise addition and multiplication, for the norm
$$
\|\hat{f}\|_{\hat{\mathscr{A}}}=\|f\|_{\mathscr{A}},\quad \forall f\in\mathscr{A}.
$$
The algebra $\mathscr{A}$ gives a general mathematical framework for the analysis of distributed delays, as it will be explained in the next subsection. Before this, let us recall the concept of Bounded Input-Bounded Output stability \cite{desoervidyasagar}.
\newdefinition{defn}{Definition}
\begin{defn}\label{def2-1-1}
A convolution system in the form~(\ref{eq2-1-1}) is said to be BIBO stable if $f\in\mathscr{A}$.
\end{defn}
More generally, we may be interested by $\mathscr{L}_{p}$-stability, for $1\leq p \leq \infty$. We say that~(\ref{eq2-1-1}) is $\mathscr{L}_{p}$-stable if for any $u(\cdot)\in\mathscr{L}_{p}(\Ree_{+})$, that is $u(\cdot)$ locally integrable over $\Ree_{+}$ and
$$
\|u\|_{\mathscr{L}_{p}}=\left(\int_{\mathds{R}_{+}}{|u(t)|^{p}\,\mathrm{d}t}\right)^{1/p}<\infty,
$$
the output $y(\cdot)$ is also in $\mathscr{L}_{p}(\Ree_{+})$. Since
\begin{equation}\label{eq2-1-5}
\|y\|_{\mathscr{L}_{p}}=\|f\ast u\|_{\mathscr{L}_{p}}\leq \|f\|_{\mathscr{A}}\|u\|_{\mathscr{L}_{p}},
\end{equation}
$\mathscr{L}_{p}$-stability is equivalent to BIBO stability~\cite{desoervidyasagar}.
%====================================================================
\subsection{Distributed delays}\label{section2-2}
Let $\Iee_{a,b}=[a,b]$ be the bounded closed interval in $\Ree_{+}$, for some reals $a$ and $b$, $0 \leq a<b$. Notations $\Iee_{0,\infty}$ or $\Ree_+$ stand for $[0,\infty[$. We define $\mathscr{K}(\Iee_{a,b})$ as the set of complex valued functions $g(\cdot)$ of the form
\begin{equation}\label{eq2-2-1}
g(t) = \left\{
\begin{array}{ll}
g_{\mathds{I}_{a,b}}(t), & t\in \Iee_{a,b}\\
0, & \text{elsewhere}
\end{array}\right.
\end{equation}
where
\begin{equation}\label{eq2-2-2}
g_{\mathds{I}_{a,b}}(t)=\sum_{i\geq 0}{\sum_{j\geq 0}{c_{ij}\,t^{j}\,\mathrm{e}^{\lambda_{i} t}}},
\end{equation}
for some $c_{ij}$ and $\lambda_{i}$ in $\Cee$, and the sums are finite. In other words, $g_{\mathds{I}_{a,b}}$ is a finite linear combination of exponential-polynomials type functions, and it is in particular a continuous function. For any real valued function in $\mathscr{K}(\Iee_{a,b})$, if some $\lambda_{i}\in\Cee$ appears in the sum, then so does its conjugate $\bar{\lambda}_{i}$, and the associated coefficients $c_{ij}$ are complex conjugates. Hence, any real valued function in $\mathscr{K}(\Iee_{a,b})$ is a function generated by real linear combinations of  $t^{j}\mathrm{e}^{\sigma_{i} t}$, $t^{j}\mathrm{e}^{\sigma_{i} t}\mathrm{sin}(\beta_{k} t)$ and $t^{j}\mathrm{e}^{\sigma_{i} t}\mathrm{cos}(\beta_{k} t)$, for some real numbers $\sigma_i$, $\beta_k$, the sums being finite. The formal definition of distributed delay is made below.
\begin{defn}\label{def2-2-1}
A distributed delay is a causal convolution system, with kernel $f$ in $\mathscr{K}(\Iee_{\vartheta_1,\vartheta_2})$, for some bounded real numbers $0\leq \vartheta_1<\vartheta_2$.
\end{defn}
In other words, a distributed delay can be written like an input-output convolution operator of the form
\begin{equation}\label{eq2-2-3}
y(t) = (f\ast u)(t)=\int_{\vartheta_{1}}^{\vartheta_{2}}{f_{\mathds{I}_{\vartheta_1,\vartheta_2}}(\tau)u(t-\tau)\,\text{d}\tau},
\end{equation}
with notations introduced in~(\ref{eq2-2-2}). The set of distributed delays, denoted by $\mathscr{G}$, is a ring. Obviously, for real valued signals, the kernel $f(\cdot)$ will be a real valued function in $\mathscr{K}(\Iee_{\vartheta_1,\vartheta_2})$. In the previous definition, we restrict ourselves to define a distributed delay like a convolution operator with kernel in $\mathscr{K}(\Iee_{\vartheta_1,\vartheta_2})$. This restriction is not conservative as we will show in Section~\ref{section3-2}. Actually, all distributed delays which appeared in the literature are particular cases of this definition. This definition is based on a rational construction, as this appears explicitly using Laplace transforms. Any distributed delay $\mathscr{G}$ admits a Laplace transform, corresponding to the finite Laplace transform of its kernel $f\in\mathscr{K}(\Iee_{\vartheta_1,\vartheta_2})$, 
\begin{equation}\label{eq2-2-4}
\hat{y}(s) = \hat{f}(s) \hat{u}(s),\quad \hat{f}(s)= \int_{\vartheta_1}^{\vartheta_2}{f_{\mathds{I}_{\vartheta_1,\vartheta_2}}(\tau)\,\mathrm{e}^{-s\tau}\,\mathrm{d}\tau},
\end{equation}
where $\hat{f}\in\hat{\mathscr{G}}$ is an entire function, \emph{i.e.} holomorphic on the whole complex plane. The notion of elementary distributed delay will greatly simplify the approximation problem. Let us define the complex valued function $\theta_\lambda(\cdot)\in\mathscr{K}(\Iee_{0,\vartheta})$, for some $\lambda\in\Cee$ and $\vartheta>0$, by
\begin{equation}\label{eq2-2-5}
\theta_\lambda(t)=\left\{\begin{array}{ll}
\mathrm{e}^{\lambda\,t} \,,&  \; t\in[0,\vartheta]\\
0 						   \,,& \text{elsewhere} 
\end{array}\right.
\end{equation}
and its Laplace transform 
\begin{equation}\label{eq2-2-6}
\hat{\theta}_\lambda(s) = \frac{1-\mathrm{e}^{-(s-\lambda)\vartheta}}{s-\lambda},
\end{equation}
which is an entire function even in $s=\lambda$ where $\hat{\theta}_\lambda(\lambda)=\vartheta$. In other words, $\lambda$ is a removable singularity, and consequently $\hat{\theta}_\lambda(s)$ has no pole.
The distributed delay whose kernel is $\theta_\lambda$ is called an elementary distributed delay. The $k$th derivative $\hat{\theta}_\lambda^{(k)}(s)$ of $\hat{\theta}_\lambda(s)$ yields
\begin{equation}\label{eq2-2-7}
\hat{\theta}_\lambda^{(k)}(s) = \int_{0}^{\vartheta}{(-\tau)^k\mathrm{e}^{-(s-\lambda)\tau}\,\mathrm{d}\tau},
\end{equation}
which is still in $\hat{\mathscr{G}}$, and corresponds to the Laplace transform of the function $\theta_\lambda^k(t)=(-t)^k\mathrm{e}^{\lambda t}$ for $t\in[0,\vartheta]$, and $0$ elsewhere. From previous definitions, we can state the following lemma, which also appeared in \cite{bretheloiseauFSA} using 2D-polynomials, and which will play a central role for approximation.
\newtheorem{lem}{Lemma}
\begin{lem}\label{def2-2-2}
Any element in $\hat{\mathscr{G}}$ can be decomposed into a finite sum of Laplace transform of elementary distributed delays and its successive derivatives.
\end{lem}
\newproof{pf}{Proof}
\begin{pf}
Take any element in $\mathscr{G}$. Its kernel $g(\cdot)$ lies in $\mathscr{K}(\Iee_{\vartheta_1,\vartheta_2})$, and writes as in~(\ref{eq2-2-2}). By time translation corresponding to the lumped delay $\vartheta_1$, it is readily a linear finite combination of elementary distributed delays $\theta_\lambda(t)$ and of the functions $\theta_\lambda^k(t)$, as defined in~(\ref{eq2-2-5}) and~(\ref{eq2-2-7}), with $\vartheta=\vartheta_2-\vartheta_1$.\hfill{$\Box$}
\end{pf}
In other terms, previous result tells us that for any $\hat{g}\in\hat{\mathscr{G}}$, there exist complex polynomials $\hat{g}_{ik}\in\Cee[\mathrm{e}^{-\vartheta s}]$ with respect to the variable $\mathrm{e}^{-\vartheta s}$ and $\lambda_i\in\Cee$, in finite number, such that 
\begin{equation}\label{eq2-2-8}
\hat{g}(s)=\sum_{i,k}{\hat{g}_{ik}(\mathrm{e}^{-\vartheta s})\hat{\theta}_{\lambda_i}^{(k)}(s)},
\end{equation}
where successive derivatives are iteratively computed by
\begin{equation}\label{eq2-2-9}
\hat{\theta}_{\lambda_i}^{(k)}(s) = (-1)^k k!\frac{1-\mathrm{e}^{-(s-\lambda_i)\vartheta}-\sum_{n=1}^{k}{\frac{\vartheta^n}{n!}\mathrm{e}^{-(s-\lambda_i)\vartheta}(s-\lambda_i)^n}}{(s-\lambda_i)^{k+1}},
\end{equation}
and $\hat{\theta}_{\lambda_i}^{(k)}(s)$ are still entire functions, since
$\hat{\theta}_{\lambda_i}^{(k)}(\lambda_i)=\frac{(-1)^k\vartheta^{k+1}}{k+1}$, for any $k\geq 0$. From~(\ref{eq2-2-8}) and taking into account that practical distributed delay is a real valued operator, we see that, grouping the terms via least common multiple, any element $\hat{g}\in\hat{\mathscr{G}}$ can be put into a fraction
\begin{equation}\label{eq2-2-10}
\hat{g}(s) = \frac{n(s,\retard)}{d(s)},
\end{equation}
where $n(s,\mathrm{e}^{-\vartheta s})\in\Ree[s,\mathrm{e}^{-\vartheta s}]$ is a real quasipolynomial with respect to the algebraically independent variables $s$ and $\retard$, and $d(s)\in\Ree[s]$. Any element in the right hand side of~(\ref{eq2-2-8}) being an entire function, $\hat{g}(s)$ is also an entire function. Hence, any zero of $d(s)$ is also a zero of $n(s,\retard)$. Furthermore, the degree with respect to $s$ of $\hat{\theta}_\lambda^{(k)}(s)$ in~(\ref{eq2-2-6}) and~(\ref{eq2-2-9}), for any $k$, being strictly negative, the degree with respect to $s$ of $\hat{g}(s)$ satisfies $\mathrm{deg}_s\,n < \mathrm{deg}_s\,d$, so that its Laplace transform is strictly proper. As in \cite{bretheloiseauFruit}, the ring $\hat{\mathscr{G}}$ is the ring of those Laplace transforms of distributed delays that are rational in the variable $s$ and $\retard$, which are entire and strictly proper (with respect to $s$). This result comes from the assumption on rational kernels in~(\ref{eq2-2-2}). The kernel of any element in $\mathscr{G}$ is obviously in $\mathscr{L}_1(\Ree_+)$, so that from~(\ref{eq2-1-5}) and the definition of $\mathscr{A}$, any distributed delay is BIBO stable. Convolution of two kernels with finite support yields another kernel with finite support, and since exponential-polynomials  type functions are closed under classical product, we have in fact that $\mathscr{G}$ is a normed subalgebra of $\mathscr{A}$ for the $\|\cdot\|_\mathscr{A}$-norm. From these definitions and properties, we are now able to formulate the approximation problem of distributed delays.
%====================================================================
\section{Approximation of distributed delays}\label{section3}
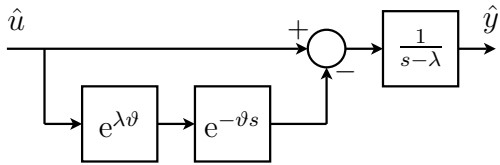
\begin{figure}[h]
\ifx\JPicScale\undefined\def\JPicScale{1}\fi
\psset{unit=\JPicScale mm}
\psset{linewidth=0.3,dotsep=1,hatchwidth=0.3,hatchsep=1.5,shadowsize=1,dimen=middle}
\psset{dotsize=0.7 2.5,dotscale=1 1,fillcolor=black}
\psset{arrowsize=1 2,arrowlength=1,arrowinset=0.25,tbarsize=0.7 5,bracketlength=0.15,rbracketlength=0.15}
\begin{pspicture}(-36,0)(65,20)
\pspolygon[](10,10)(20,10)(20,0)(10,0)
\pspolygon[](25,10)(35,10)(35,0)(25,0)
\psline{->}(0,15)(40,15)
\psline(5,15)(5,5)
\psline{->}(5,5)(10,5)
\psline{->}(20,5)(25,5)
\psline(35,5)(42.5,5)
\rput{0}(42.5,15){\psellipse[](0,0)(2.5,2.5)}
\psline{<-}(42.5,12.5)(42.5,5)
\psline{->}(45,15)(50,15)
\pspolygon[](50,20)(60,20)(60,10)(50,10)
\rput(30,5){$\mathrm{e}^{-\vartheta s}$}
\rput(15,5){$\mathrm{e}^{\lambda\vartheta}$}
\rput(55,15){$\frac{1}{s-\lambda}$}
\psline{->}(60,15)(65,15)
\rput(1.25,19){$\hat{u}$}
\rput(64.38,19){$\hat{y}$}
\rput(38.75,17.5){\small $+$}
\rput(45,11.88){\small $-$}
\end{pspicture}
\caption{Realization of an element $\hat{\theta}_\lambda(s)\in\mathscr{K}_s(\Iee_{0,\vartheta})$, with $\mathrm{Re}\,\lambda<0$.}\label{fig3-1}
\end{figure}

For this purpose, we define two subspaces in direct sum in $\mathscr{K}(\Iee_{\vartheta_1,\vartheta_2})$, denoted respectively $\mathscr{K}_s(\Iee_{\vartheta_1,\vartheta_2})$ and $\mathscr{K}_u(\Iee_{\vartheta_1,\vartheta_2})$, consisting of linear combinations of exponential-polynomials type functions on some finite interval as in~(\ref{eq2-2-2}), with $\mathrm{Re}\,\lambda_i<0$ and $\mathrm{Re}\,\lambda_i\geq 0$, respectively, for all $i\geq 0$.\\ 
From Lemma~\ref{def2-2-2} and~(\ref{eq2-2-8}), we know that any element $g\in\mathscr{K}_s(\Iee_{\vartheta_1,\vartheta_2})$ is a linear combination of elements $\theta_{\lambda_i}^{k}$ in $\mathscr{K}_s(\Iee_{\vartheta_1,\vartheta_2})$, for which $\mathrm{Re}\,\lambda_i<0$. Hence, any element in $\mathscr{K}_s(\Iee_{\vartheta_1,\vartheta_2})$ can be numerically implemented with stability using elementary blocks, as illustrated in Fig.~{\ref{fig3-1}} for $\theta_\lambda(t)$. Note that in practice, since $\lambda\in\Cee$, we should implement it using a real decomposition.
Approximation for the subclass $\mathscr{K}_s(\Iee_{\vartheta_1,\vartheta_2})$ of distributed delays is not required, since this implementation is realized with stability, and requires, from an implementation point of view, only two pointwise delays, namely $0$ and $\vartheta$. However, this is no more true for elements in $\mathscr{K}_u(\Iee_{\vartheta_1,\vartheta_2})$, where numerical realization for $\theta_\lambda(t)$ as in Fig.~{\ref{fig3-1}} yields to an unstable system. Note that another realization for $\hat{\theta}_\lambda(s)$ may be
$$
\dot{x}(t) = \lambda x(t) + u(t)-\mathrm{e}^{\lambda\vartheta}u(t-\vartheta),
$$
which is however still numerically unstable for $\mathrm{Re}\,\lambda\geq 0$. Intuitively, this instability is a consequence of a non exact numerical cancellation of the removable singularity $s=\lambda$ of $\hat{\theta}_\lambda(s)$. Therefore, we will focus our attention on distributed delays whose kernels are in $\mathscr{K}_u(\Iee_{\vartheta_1,\vartheta_2})$.

\noindent The question can now be stated as follows: How to define a continuous time approximation of a distributed delay whose kernel lies in $\mathscr{K}_u(\Iee_{\vartheta_1,\vartheta_2})$? We give complete answers to this question in the next two subsections. We start by some general remarks on approximation over $\mathscr{A}$, and hence give answers to the previous question.
%=====================================================================
\subsection{General comments on approximation}\label{section3-1}
\begin{figure}[h!]
\begin{center}
\includegraphics[width=3.5in,height=2.5in]{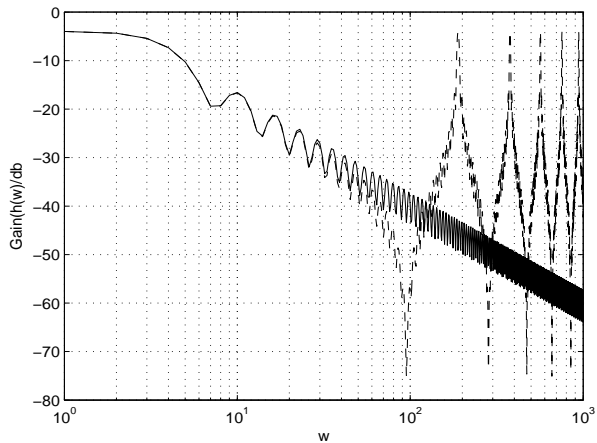}
\end{center}
\caption{Bode gain diagram of $\hat{\theta}_1(s)$ and an approximation based on integral approximation (trapezoidal Newton-Cotes method with 30 nodes).}\label{fig3-1-1}
\end{figure}

Any distributed delay with kernel $f\in\mathscr{K}(\Iee_{\vartheta_1,\vartheta_2})$ is a convolution operator. Its $\mathscr{L}_\infty$-induced norm is
$$
\|f\|_\mathscr{A} = \|f\|_{\mathscr{L}_1}.
$$
From~(\ref{eq2-1-5}), we see that this norm is actually an upper bound of all its induced $\mathscr{L}_p$-norms, for $1\leq p\leq\infty$. A natural metric is then obtained from the $\mathscr{A}$-norm, which in turn can be enclosed in the graph topology. Consider a distribution $p$ with a coprime factorization $(n,d)$ in $\mathcal{A}$, that is $\hat{p}=\hat{n}\hat{d}^{-1}$, and such that there exist $x$ and $y$ in $\mathscr{A}$ satisfying 
\begin{equation}\label{eq3-1-1}
n\ast x + d\ast y = \delta,
\end{equation}
or in the Laplace domain
\begin{equation}\label{eq3-1-2}
\hat{n}\hat{x}+\hat{d}\hat{y}=1,\;\forall s\in\Cee,\; \mathrm{Re}\,s\geq 0.
\end{equation}
A neighbourhood of $p$ in the graph topology is the set of all plants of the form $n_\Delta d_\Delta^{-1}$ where $(n_\Delta,d_\Delta)$ belongs to some ball in $\mathscr{A}$ centered at $(n,d)$. If $p$ is itself in $\mathscr{A}$, for any positive number $r$, a neighbourhood of $p$ in the graph topology writes
$$
\mathscr{B}(p,r) = \{p_\Delta\in\mathscr{A}: \; \|p-p_\Delta\|_\mathscr{A}\leq r\}.
$$
In~\cite[Ch. 7]{vidyasagarbook}, it was shown that the graph topology is the weakest topology in which feedback stability is a robust property. Graph topology being metrizable and hence first-countable, for BIBO-stable systems, norm topology and graph topology are the same. Let $\varepsilon$ be a given suitable small positive number.  Hence, we say that a distribution $p_\varepsilon$ in $\mathcal{A}$ is an approximation of $p$ in the graph topology if the element $p_\varepsilon$ is close to $p$, that is $p_\varepsilon\in\mathscr{B}(p,\varepsilon)$, or equivalently $\|p-p_\varepsilon\|_\mathscr{A}\leq\varepsilon$. This property can be related to convergence in $\mathscr{A}$, but we need first to define the class of operators that will approximate $p$ in the graph topology. The most commonly used class is the so-called class of lumped systems, which consists in convolution operators whose kernels lie in the set $\mathscr{Q}$, defined by  
\begin{equation}\label{eq3-1-3}
\mathscr{Q} = \left\{g\in\mathscr{A}: \; g(t)=g_0\delta(t)+g_a(t),\;g_0\in\Cee,\;g_a\in\mathscr{K}_s(\Iee_{0,\infty})\right\}.
\end{equation}
From~\cite{vidyasagaranderson}, we know that any plant $p$ with coprime factorization $(n,d)$ over $\mathscr{A}$ can be approximated by a lumped plant in the graph topology if and only if there exists a real constant matrix $M$ of rank $1$ such that
\begin{equation}\label{eq3-1-4}
MS_{pa}=0,\; \text{where} \; S_{pa}=\left[\begin{matrix}n_{pa}\\d_{pa}\end{matrix}\right],
\end{equation}
and $(\cdot)_{pa}$ denotes the purely atomic part. Distributed delays having strictly proper transfer functions, we are interested for approximation with the class of strictly proper lumped systems, denoted $\mathscr{Q}_{s}$, set of elements in  $\mathscr{K}_s(\Iee_{0,\infty})$. Such elements satisfy ordinary differential equations, and are easily simulated. From \cite{kammler}, we know that the closure of $\mathscr{Q}_s$ is $\mathscr{L}_1(\Ree_+)$, since every function in $\mathscr{L}_1(\Ree_+)$ can be approximated by a sum of exponentials. This yields the following result, that can be seen as a particular case of the above given condition~(\ref{eq3-1-4}).
\newtheorem{thm}{Theorem}
\begin{thm}\label{def3-1-1}
Any plant in $\mathscr{L}_1(\Ree_+)$ can be approximated by a lumped system in the graph topology.
\end{thm}
\begin{pf}
Let $f$ be an element in $\mathscr{L}_1(\Ree_+)$. A coprime factorization of $f$ over $\mathscr{A}$ is $n=f$ and $d=\delta$. Since $n_{pa}=0$, any real constant matrix $M=\left[0 \;\; \star\right]$, $(\star)$ denoting an arbitrary real number, satisfies~(\ref{eq3-1-4}). Hence, from~\cite[Th. 4.1]{vidyasagaranderson}, $f$ can be approximated by a lumped system in the graph topology.\hfill{$\Box$}
\end{pf}
Note that if the element $f$ in $\mathscr{L}_1(\Ree_+)$ is strictly proper, then an approximation in the graph topology realized by a lumped system will be in fact in the set $\mathscr{Q}_s$. Any convolution system whose kernel lies in $\mathscr{L}_1(\Ree_+)$ can be approximated by a lumped system. Roughly speaking, since any distributed delay is a strictly proper fraction and is BIBO-stable, we know that it can be approximated in the graph topology by a strictly proper lumped system.
\newtheorem{cor}{Corollary}
\begin{cor}\label{def3-1-2}
Any distributed delay can be approximated by a lumped system in the graph topology.
\end{cor}
Others approximation classes can be used to numerically realize a distributed delay. In the literature, appears a method based on a numerical integral approximation, but such a method was shown to be at the core of various undesired results, like instability or numerical sensitivity. Let us briefly make some considerations on such approximations, introducing what we called input-output approximation, and let us show why such approximations are not suitable for our problem. For any distributed delay $y(t)=(f\ast u)(t)$, we say that $y_{\mathrm{app}}$ is an $\mathscr{L}_p$ input-output approximation of $y$, if for some arbitrarily $\varepsilon>0$, $\left\|y-y_{\mathrm{app}}\right\|_{\mathscr{L}_{p}}\leq\varepsilon$ holds. For instance, consider the $\mathscr{L}_\infty$ input-output approximation based on integral approximation, like Newton-Cotes quadrature methods, which leads to an approximation $y_{\mathrm{app}}(t)$ of the form
\begin{equation}\label{eq3-1-5}
y_{\mathrm{app}}(t) = \sum^{q}_{k=0}c_{k}f(\tau_{k})u(t-\tau_{k}),
\end{equation}
where reals $c_{k}$ and $\tau_k$ depend on the applied method~\cite[Ch. 4]{ralstonrabinowitz}. Various problems arise with this kind of approximation. First note that such an approximation is made for a given $u(\cdot)$. Hence, if such an input changes, the properties of this approximation failed, in general. Note also that, as mentioned by~\cite{mirkinApproximation}, the Laplace transform of~(\ref{eq3-1-5}) yields to an input-output approximation transfer function  that is proper, since it contains only pointwise delays. From~(\ref{eq2-2-8}) and~(\ref{eq2-2-10}), a distributed having a strictly proper transfer function, this approximation has a poor accuracy for high frequencies, which will have as consequence a high sensitivity to high-frequency plant uncertainties, in particular for the closed-loop system. This fact is illustrated in Fig.~\ref{fig3-1-1}, where a frequency diagram of $\hat{\theta}_\lambda(s)$
and an input-output approximation is plotted. This negative result is directly interpretable in the graph topology by the following result.

\begin{thm}\label{def3-1-3}
Any distributed delay can not be approximated in the graph topology by a purely atomic distribution.
\end{thm}
\begin{pf}
Consider an element in $\mathscr{G}$ with kernel $g$. Then $g\in\mathscr{L}_1(\Ree_+)$. Consider an arbitrary purely atomic distribution with pointwise delays of the form
$$
p(t)=\sum_{i\geq 0}{p_i\delta(t-t_i)}
$$
where $\sum_{i\geq 0}{|p_i|}<\infty$. Such a distribution corresponds to the kernel of the approximation using numerical integral approximation like in~(\ref{eq3-1-5}). Then
$$
\|g-p\|_\mathscr{A} = \|g\|_{\mathscr{L}_1}+\sum_{i\geq 0}{|p_i|}.
$$
Hence $\|g-p\|_\mathscr{A} \geq \|g\|_{\mathscr{L}_1}$, and consequently $p$ can not approximate $g$ in the graph topology.
\hfill{$\Box$}
\end{pf}
It is clear that input-output approximation is not suitable for numerical implementation. This was illustrated in~\cite{vanassche},~\cite{zhongphd}, or~\cite{mirkinApproximation}. Hence we focus our attention to the notion of approximation defined in the graph topology. This approximation is equivalent to realize an approximation of the kernel of the convolution, so we may also call it kernel approximation, to discern with respect to input-output approximation. 
%====================================================================
\subsection{Kernel approximations}\label{section3-2}
Approximation in the graph topology for an element $f$ in $\mathscr{K}_s(\Iee_{\vartheta_1,\vartheta_2})$ requires to find, for some given $\varepsilon>0$, an element $f_\mathrm{app}$, defined in general in $\mathscr{A}$, such that $f_\mathrm{app}\in\mathscr{B}(f,\varepsilon)$.
Element $f$ being in particular in $\mathscr{L}_1(\Ree_+)$, adding purely atomic part in approximation can not improve approximation, so we realize it over the Banach algebra $\mathscr{L}_1(\Ree_+)$. Indeed, if an approximation writes
$$
f_{\mathrm{app}} = f_a^{\mathrm{app}}+f_{pa}^{\mathrm{app}},
$$
where $f_a^{\mathrm{app}}$ and $f_{pa}^{\mathrm{app}}$ stand for the atomic and purely atomic parts of $f_{\mathrm{app}}$, respectively, the following decomposition holds
$$
\|f-f_\mathrm{app}\|_\mathscr{A} = \|f-f_a^{\mathrm{app}}\|_{\mathscr{A}}+\|f_{pa}^{\mathrm{app}}\|_{\mathscr{A}}.
$$
The purely atomic part of the approximation is in an independent sum, and can then be reduced to zero to reduce the approximation error. We have seen that an approximation can be obtained by lumped systems, with rational and stable transfer functions. Since we want to obtain an approximation which is global with respect to the time, we are interested in uniform convergent sequence of functions in $\mathscr{L}_1(\Ree_+)$ to $f$.  Consider a sequence $f_{n}$ of functions which uniformly converges to $f$, that is, for any $\varepsilon>0$, there exists $n$ in $\Nee$ such that for any $\kappa\geq n$, $\|f-f_{\kappa}\|_\mathscr{A}\leq\varepsilon$. So, for any $\kappa\geq n$, $f_{\kappa}$ is a kernel approximation of $f$ in the graph topology. We propose next a method to approximate an element in $\mathscr{K}_u(\Iee_{\vartheta_1,\vartheta_2})$ by elements in $\mathscr{K}_s(\Iee_{\vartheta_1,\vartheta_2})$. Since these last elements are easily implementable, this will give us an effective approximation. In other words, we approximate a distributed delay by another distributed delays. We start by considering the case of an elementary distributed delay. 
\begin{lem}\label{def3-2-1}
Any distributed delay with kernel $\theta_\lambda(\cdot)$ in $\mathscr{K}_u(\Iee_{0,\vartheta})$ can be approximated by distributed delays with kernels in $\mathscr{K}_s(\Iee_{0,\vartheta})$ for the graph topology.
\end{lem}
\begin{pf}
Let $\theta_\lambda(\cdot)$ be the kernel of an elementary distributed delay, with $\mathrm{Re}\,\lambda\geq 0$. Let $\mu\in\Cee$ such that $\mathrm{Re}\,\mu=\alpha<0$. We define the transform
\begin{equation}\label{eq3-2-1}
\Theta_\lambda(s)=(\alpha s)^{-1}\,\theta_\lambda(-\alpha^{-1}\,\mathrm{ln}\,s),\quad \mathrm{e}^{-\vartheta\alpha}\leq s\leq 1.
\end{equation} 
The function $\Theta_\lambda(\cdot)$ is continuous over $[\mathrm{e}^{-\vartheta\alpha},1]$, and we have
\begin{equation}\label{eq3-2-2}
\|\Theta_\lambda\|_{\mathscr{L}_1} = \|\theta_\lambda\|_\mathscr{A}.
\end{equation}
By M\"untz-Sz\'asz theorem, we see that $\Theta_\lambda(\cdot)$ can be approximated with respect to $\mathscr{L}_1$-norm, as closed as desired, by a function of the form
\begin{equation}\label{eq3-2-3}
\Theta_{\lambda,\mathrm{app}}(s) = (\alpha s)^{-1}\, \Psi_{\lambda,\mathrm{app}}(s),
\end{equation}
where $\Psi_{\lambda,\mathrm{app}}$ is a polynomial in $s$. Since norms are preserved in~(\ref{eq3-2-2}), this comes down to approximate $\theta_\lambda$ as closely as we please, using inverse transform of~(\ref{eq3-2-1}), by a sum of exponentials of the form
$$
\psi_{\lambda,\mathrm{app}}(t) = \Psi_{\lambda,\mathrm{app}}(\mathrm{e}^{-\alpha t}),
$$
which is clearly an element in $\mathscr{K}_s(\Iee_{0,\vartheta})$.
\hfill{$\Box$}
\end{pf}
In other words, an approximation of $\theta_\lambda$ in $\mathscr{A}$ can be a function of the form
\begin{equation}\label{eq3-2-4}
\psi_{\lambda,n}(t) = \sum_{i=1}^{n}{\gamma_{i,n}\theta_{\alpha_i}(t)},
\end{equation}
with $\alpha_i\in\Cee$, $\mathrm{Re}\,\alpha_i<0$ some arbitrarily complex numbers, and $\gamma_{i,n}$ some suitable constants. The order of the approximation $n$ describes the number of parallel distributed elements to be added, to get the desired approximation accuracy. Such an approximation does not increase the degree of the elementary distributed delay, since it is a sum of distributed delays which are strictly proper and have same degree than $\hat{\theta}_\lambda(s)$. From the synthesis of an approximation for the elementary distributed delay, we show in the next two results that it allows to construct an explicit approximation for any distributed delay.
\begin{lem}\label{def3-2-2}
Let $\psi_{\lambda,n}$ be an element in $\mathscr{K}_s(\Iee_{0,\vartheta})\cap\mathscr{B}(\theta_\lambda,\varepsilon)$, for a given $\varepsilon>0$, and $k$ in $\Nee$. Then
$\psi_{\lambda,n}^k\in\mathscr{K}_s(\Iee_{0,\vartheta})\cap\mathscr{B}(\theta_\lambda^k,\vartheta^k\varepsilon)$, where
$$
\psi_{\lambda,n}^k(t)=\left\{\begin{array}{ll}
(-t)^k\psi_{\lambda,n}(t) \,,&  \; t\in[0,\vartheta]\\
0 						   \,,& \text{elsewhere} 
\end{array}\right..
$$
\end{lem}
\begin{pf}
Let $\varepsilon>0$, and take $\psi_{\lambda,n}$ an element in $\mathscr{K}_s(\Iee_{0,\vartheta})\cap\mathscr{B}(\theta_\lambda,\varepsilon)$, that is 
$$
\psi_{\lambda,n}(t) = \sum_{i=1}^{n}{\gamma_{i,n}\,\theta_{\alpha_i}(t)},
$$
with $\gamma_{i,n}$ in $\Cee$, and $\|\theta_\lambda-\psi_{\lambda,n}\|_\mathscr{A}\leq\varepsilon$. By Laplace transform and $k$th order differentiation, we get
$$
\hat{\psi}_{\lambda,n}^{(k)}(s) = \sum_{i=1}^{n}{\gamma_{i,n}\,\hat{\theta}_{\alpha_i}^{(k)}(s)}.
$$
In the time domain, this last identity corresponds to
$$
\psi^k_{\lambda,n}(t) = \sum_{i=1}^{n}{\gamma_{i,n}\,\theta^k_{\alpha_i}(t)}=(-t)^k\psi_{\lambda,n}(t).
$$
Clearly $\psi_{\lambda,n}^k\in\mathscr{K}_s(\Iee_{0,\vartheta})$. Hence, it remains to show that $\psi_{\lambda,n}^k\in\mathscr{B}(\theta_\lambda^k,\vartheta^k\varepsilon)$. But  
$\psi_{\lambda,n}\in\mathscr{B}(\theta_\lambda,\varepsilon)$, so we have
\begin{eqnarray*}
\|\theta_\lambda^k-\psi_{\lambda,n}^k\|_\mathscr{A} & = & \int_{0}^{\vartheta}{\left|(-t)^k(\theta_\lambda(t)-\psi_{\lambda,n}(t))\right|\,\text{d}t}\\
& \leq &   \vartheta^k\|\theta_\lambda-\psi_{\lambda,n}\|_\mathscr{A}\leq \vartheta^k \varepsilon,
\end{eqnarray*}
which completes the proof.
\hfill{$\Box$}
\end{pf}
For any distributed delay, we have the following general result.
\begin{thm}\label{def3-2-3}
Let $\hat{g}\in\hat{\mathscr{G}}$ be an arbitrary distributed delay, of the form
$$
\hat{g}(s) = \sum_{i,k}{\hat{g}_{ik}(\retard)\hat{\theta}_{\lambda_i}^{(k)}(s)},
$$
and let $\varepsilon_i$ be given positive real numbers. For any elements $\psi_{\lambda_i,n}$ in 
$\mathscr{K}_s(\Iee_{0,\vartheta})\cap\mathscr{B}(\theta_{\lambda_i},\varepsilon_i)$, we define
$$
\hat{g}_{\mathrm{app}}(s) = \sum_{i,k}{\hat{g}_{ik}(\retard)\hat{\psi}_{\lambda_i,n}^{(k)}(s)}.
$$
Then $g_\mathrm{app}$ lies in $\mathscr{K}_s(\Iee_{0,\vartheta})\cap\mathscr{B}(g,\tilde{\varepsilon})$, where $\tilde{\varepsilon}=M \cdot \mathop{\mathrm{max}}_i\,\varepsilon_i$ for some positive constant $M$.
\end{thm}
\begin{pf}
Let $g$ be the kernel of a distributed delay. According to~(\ref{eq2-2-8}), we decompose it as a linear combination of elementary distributed delays and their successive derivatives. From the definition of $g_{\mathrm{app}}$, which is clearly in $\mathscr{G}$, we have 
\begin{eqnarray*}
\|g-g_{\mathrm{app}}\|_\mathscr{A} & = & \left\|\sum_{i,k}{\hat{g}_{ik}(\retard)(\hat{\theta}_{\lambda_i}^{(k)}(s)-\hat{\psi}_{\lambda_i,n}^{(k)}(s))}\right\|_{\hat{\mathscr{A}}}\\
& \leq & \sum_{i,k}{\|\hat{g}_{ik}(\retard)\|_{\hat{\mathscr{A}}}
\|\hat{\theta}_{\lambda_i}^{(k)}(s)-\hat{\psi}_{\lambda_i,n}^{(k)}(s)}\|_{\hat{\mathscr{A}}}.
\end{eqnarray*}
Since, for any $i$, $\psi_{\lambda_i,n}$ are in 
$\mathscr{B}(\theta_{\lambda_i},\varepsilon_i)$, using Lemma~\ref{def3-2-2}, we get 
\begin{eqnarray*}
\|g-g_{\mathrm{app}}\|_\mathscr{A} & \leq & \sum_{i,k}{\vartheta^k \,\|\hat{g}_{ik}(\retard)\|_{\hat{\mathscr{A}}}\,\|\hat{\theta}_{\lambda_i}(s)-\hat{\psi}_{\lambda_i,n}(s)}\|_{\hat{\mathscr{A}}}\\
& \leq & \sum_{i,k}{\vartheta^k\varepsilon_i\,\|\hat{g}_{ik}(\retard)\|_{\hat{\mathscr{A}}}}.
\end{eqnarray*}
Denote the positive bounded constant $M = \sum_{i,k}{\vartheta^k\|\hat{g}_{ik}(\retard)\|_{\hat{\mathscr{A}}}}$.
Hence
$$
\|g-g_{\mathrm{app}}\|_\mathscr{A} \leq \tilde{\varepsilon} = M\cdot \mathop{\mathrm{max}}_i\,\varepsilon_i,
$$ 
that is $g_{\mathrm{app}}$ lies in $\mathscr{B}(g,\tilde{\varepsilon})$.
\hfill{$\Box$}
\end{pf}
By suitable choices for $\varepsilon_i$, the upper bound $\tilde{\varepsilon}$ can be reduced arbitrarily, so that we can find $g_{\mathrm{app}}$ as close as we please of $g$. Note that we particularize the proofs of previous results with approximations with kernels in $\mathscr{K}_s(\Iee_{0,\vartheta})$, but they can be trivially extended to others approximations over $\mathscr{A}$. Previous results state that from the approximation of the elementary distributed delay, we can realize a kernel approximation in the graph topology for any distributed delay in $\mathscr{G}$.  Previous results state that, for any $f\in\mathscr{K}(\Iee_{\vartheta_1,\vartheta_2})$ and $\varepsilon>0$, there exists an element $f_{\text{app}}\in\mathscr{K}_s(\Iee_{\vartheta_1,\vartheta_2})$ such that
\begin{equation}\label{eq3-2-5}
\|f-f_{\text{app}}\|_\mathscr{A}\leq\varepsilon.
\end{equation}
Said differently,  we have the following corollary.
\begin{cor}\label{def3-2-4}
The set $\mathscr{K}_s(\Iee_{\vartheta_1,\vartheta_2})$ is dense in $\mathscr{K}(\Iee_{\vartheta_1,\vartheta_2})$ for the graph topology.
\end{cor}
\begin{pf}
Obvious from Lemma~\ref{def2-2-2} and Theorem~\ref{def3-2-3}.\hfill{$\Box$}
\end{pf}
In the proof of Lemma~\ref{def3-2-1}, convergence of polynomial approximation is uniform, so we get here a uniform convergence of this  approximation. The assumption that kernel of distributed delay lie in $\mathscr{K}(\Iee_{\vartheta_1,\vartheta_2})$ is not restrictive. Indeed, if this is not the case, any continuous function $g$ in $[\vartheta_1,\vartheta_2]$ can be approximated with respect to $\mathscr{L}_1$-norm as close as we please by a function in $\mathscr{K}(\Iee_{\vartheta_1,\vartheta_2})$~\cite{kammler}. Previous results addressed the approximation of distributed delays using kernels in $\mathscr{K}_s(\Iee_{\vartheta_1,\vartheta_2})$, which is based on uniform convergence of polynomials to any continuous function over $\Iee_{\vartheta_1,\vartheta_2}$. Other approximation can be proposed.  
Indeed, any element $f$ in $\mathscr{K}(\Iee_{\vartheta_1,\vartheta_2})$ can also be uniformly approximated by stepwise continuous functions in $\mathscr{L}_1(\Iee_{0,\infty})$. This can be seen as a particular case of the previous approximation. For this, consider
\begin{equation}\label{eq3-2-6}
\theta_{\lambda,\mathrm{app}}(t) = \sum_{i=0}^{n}{\gamma_i\,\psi(t-t_i)},
\end{equation}
where $t_i=\frac{i\vartheta}{n}$, and $\psi(\cdot)$ is a function in $\mathscr{K}_s(\Iee_{0,\infty})$. Then
\begin{eqnarray*}
\|\theta_\lambda-\theta_{\lambda,\mathrm{app}}\|_\mathscr{A} & = &
\int_{0}^{\vartheta}{|\theta_\lambda(t)-\theta_{\lambda,\mathrm{app}}(t)|\,\mathrm{d}t} + \int_{\vartheta}^{\infty}{|\theta_{\lambda,\mathrm{app}}(t)|\,\mathrm{d}t} \\
& = & \sum_{k=0}^{n-1}{\int_{k\frac{\vartheta}{n}}^{(k+1)\frac{\vartheta}{n}}{|\theta_\lambda(t)-\sum_{i=0}^{k}{\gamma_i\psi(t-t_i)}|\,\mathrm{d}t}}+\int_{\vartheta}^{\infty}{|\theta_{\lambda,\mathrm{app}}(t)|\,\mathrm{d}t}.
\end{eqnarray*}
From Lemma~\ref{def3-2-1}, we see that approximation~(\ref{eq3-2-6}) uniformly converges to $\theta_\lambda(\cdot)$ in the graph topology by a suitable choice of coefficients $\gamma_i$. Taking $\psi(t)=\mathrm{e}^{-\alpha t}h(t)$, where $\alpha>0$ and $h(\cdot)$ stands for the Heaviside function, we obtain as a particular case the result obtained in~\citep{mirkinApproximation}, and separately in~\citep{mondiemichiels}, where a low pass filter is added in the integral approximation with lumped delayed distributions. Indeed, Laplace transform of~(\ref{eq3-2-6}) is of the required form
$$
\hat{\theta}_{\lambda,\mathrm{app}}(s) = \frac{1}{s+a}\sum_{i=0}^{n}{\gamma_i\,\mathrm{e}^{-st_i}}
$$
for $a>0$. Coefficients $\gamma_i$ can be obtained for instance from some numerical integral approximation, to guarantee approximation over $\mathscr{A}$, or equivalently over $\mathscr{L}_1(\Iee_{0,\infty})$.
%====================================================================
\section{Discussion on approximation}\label{section4}
We particularize in this section a proposal for approximation, and we analyze the main properties of the error in the time and frequency domains. These properties will be still valid for any other approximation in the graph topology over $\mathscr{A}$.

\subsection{Proposal for a constructive approximation}\label{section4-1}
\begin{figure}[h!]
\begin{center}
\includegraphics[width=3.5in,height=2.5in]{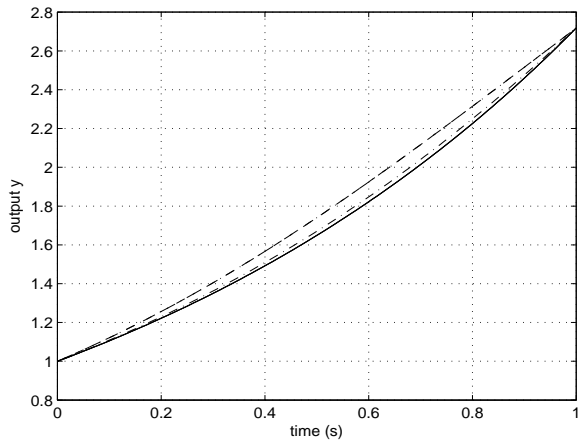}
\end{center}
\caption{Kernel approximations $\theta_{1,\mathrm{app}}(t)$ using exponentials in the time domain of the kernel $\theta_1(t)$ (continuous line), for orders $n=5$ (dashed) and $n=10$ (dot-dashed).}\label{fig4-1-1}
\end{figure}

\begin{figure}[h!]
\begin{center}
\includegraphics[width=3.5in,height=2.5in]{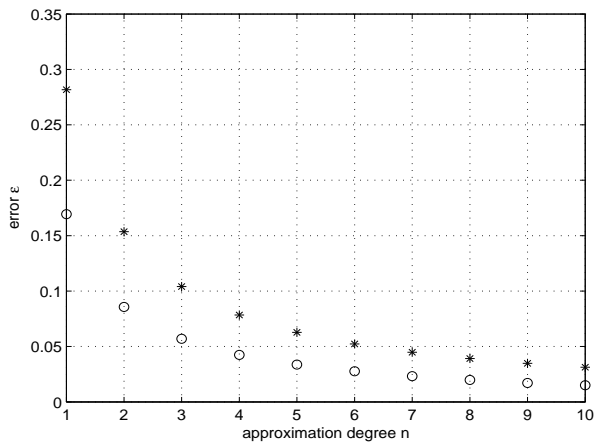}
\end{center}
\caption{Representation of the norm error $\|\theta_1-\theta_{1,n}\|_\mathscr{A}$ with respect to the order $n$ of the approximation, for $\alpha=1$ (star) and $\alpha=1/5$ (circle).}\label{fig4-1-2}
\end{figure}

We propose in this subsection an effective and constructive approximation using elements in $\mathscr{K}_s(\Iee_{\vartheta_1,\vartheta_2})$. We start this proposal by the particular case of $\theta_0(\cdot)$. We denote $\mathrm{C}^k_n = \frac{n!}{k!\,(n-k)!}$, $\psi_0(\xi)=\theta_0(-\alpha^{-1}\,\mathrm{ln}\,\xi)$, for $\xi\in]0,1]$, and $\psi_0(0)=0$. 
\begin{lem}\label{def4-1-1}
Consider the sequence $\theta_{0,n}(\cdot)$ in $\mathscr{K}_s(\Iee_{0,\vartheta})$ described by 
$$
\theta_{0,n}(t) =  \sum_{k=0}^{n}{\mathrm{C}^k_n\,\psi_0\left(\frac{k}{n}\right)\mathrm{e}^{-\alpha kt}(1-\mathrm{e}^{-\alpha t})^{n-k}},\; t\in\Iee_{0,\vartheta}
$$
and 0 elsewhere. Then $\theta_{0,n}(\cdot)$ uniformly converges to $\theta_0(\cdot)$ for the $\mathscr{A}$-norm.
\end{lem}
\begin{pf}
Define
\begin{equation}\label{eq4-1-1}
\psi_0(\varrho)=\theta_0(-\alpha^{-1}\,\mathrm{ln}\,\varrho),\quad \varrho\in]0,1],
\end{equation}
and $\psi_0(0)=0$. The function $\psi_0(\cdot)$ has a bounded step discontinuity for $\varrho=\mathrm{e}^{-\alpha\vartheta}$. Since $\mathrm{e}^{-\alpha\vartheta}$ is irrational, there exists $\kappa$ in $0,1,\ldots,n-1$ such that $\kappa/n<\mathrm{e}^{-\alpha\vartheta}<(\kappa+1)/n$. We take $\psi_{0,c}(\cdot)$ any continuous function over $[0,1]$ satisfying
$$
\psi_{0,c}(\varrho) = \psi_0(\varrho), \varrho\in[0,\kappa/n]\cup[(\kappa+1)/n,1],
$$
that will remove this discontinuity and will be close as desired to $\psi_0(\cdot)$ for the $\mathscr{L}_1(\Iee_{0,1})$-norm. In particular $\psi_{0,c}(\cdot)$ will satisfy
$$
\psi_{0,c}(k/n)=\psi_0(k/n),\; k=0, 1,\ldots,n.
$$
We now approximate as close as we please the function $\psi_{0,c}$ (or equivalently $\psi_0$) by 
Bernstein polynomials
$$
\psi_{0,\mathrm{app}}(\varrho) =  \sum_{k=0}^{n}{\mathrm{C}^k_n\,\psi_0\left(\frac{k}{n}\right)\varrho^k(1-\varrho)^{n-k}}.
$$
By the inverse transform of~(\ref{eq4-1-1}), we obtain
$$
\theta_{0,n}(t) = \psi_{0,\mathrm{app}}(\mathrm{e}^{-\alpha t}),
$$
which is in $\mathscr{K}_s(\Iee_{0,\vartheta})$ since $\psi_0(0)=0$. Furthermore, from Lemma~\ref{def3-2-1}, we know that $\theta_{0,n}(\cdot)$ converges uniformly, by construction, to $\theta_0(\cdot)$.
\hfill{$\Box$}
\end{pf}
Previous lemma gives an approximation of $\theta_0$ in $\mathscr{B}(\theta_0,\varepsilon_n)$, where positive upper bound $\varepsilon_n$ can be reduced arbitrarily when $n$ increases. This algorithm can be easily generalized to approximations in $\mathscr{K}(\Iee_{0,\infty})$ or of the form~(\ref{eq3-2-6}), by modifying the domain of definition for~(\ref{eq4-1-1}). For the more general case of $\theta_\lambda(\cdot)$, we propose the following constructive solution.
\begin{lem}\label{def4-1-2}
The sequence of functions defined over $\mathscr{K}_s(\Iee_{0,\vartheta})$ by
$$
\theta_{\lambda,n}(t) = \frac{1}{\left(1-\mathrm{e}^{-\alpha\vartheta}\right)^n}\sum_{k=0}^{n}{\mathrm{C}^k_n\,\Phi_\lambda\left(\frac{k}{n}\right)\left(\mathrm{e}^{-\alpha t}-\mathrm{e}^{-\alpha\vartheta}\right)^k\left(1-\mathrm{e}^{-\alpha t}\right)^{n-k}},
$$
for $t\in\Iee_{0,\vartheta}$, and 
$$
\Phi_\lambda(\mu) = \theta_\lambda(-\alpha^{-1}\,\mathrm{ln}((1-\mathrm{e}^{-\alpha\vartheta})\mu+\mathrm{e}^{-\alpha\vartheta})),\; \mu\in[0,1],
$$
uniformly converges to $\theta_\lambda(\cdot)$.
\end{lem}
\begin{pf}
From Lemma~\ref{def3-2-1}, define
$$
\Phi_\lambda(\mu) = \theta_\lambda(-\alpha^{-1}\,\mathrm{ln}((1-\mathrm{e}^{-\alpha\vartheta})\mu+\mathrm{e}^{-\alpha\vartheta}))
$$
with $\mu\in[0,1]$. Polynomial approximation $\Phi_{\lambda,\mathrm{app}}(\cdot)$ can be obtained from Bernstein polynomials, that is
\begin{equation}
\Phi_{\lambda,\mathrm{app}}(\mu) = \sum_{k=0}^{n}{\mathrm{C}^k_n\,\Phi_\lambda\left(\frac{k}{n}\right)\mu^k(1-\mu)^{n-k}}.
\end{equation} 
The inverse transformation for $\Phi_{\lambda,\mathrm{app}}$ yields to
$$
\theta_{\lambda,n}(t) = \frac{1}{\left(1-\mathrm{e}^{-\alpha\vartheta}\right)^n}\sum_{k=0}^{n}{\mathrm{C}^k_n\,\Phi_\lambda\left(\frac{k}{n}\right)\left(\mathrm{e}^{-\alpha t}-\mathrm{e}^{-\alpha\vartheta}\right)^k\left(1-\mathrm{e}^{-\alpha t}\right)^{n-k}},
$$
which uniformly converges to $\theta_\lambda(\cdot)$.
\hfill{$\Box$}
\end{pf}
This approximation writes like a sum of elements in $\mathscr{K}_s(\Iee_{0,\vartheta})$ and $\theta_0(\cdot)$. With Lemma~\ref{def4-1-1} and approximation of $\theta_0(\cdot)$, such an approximation is defined over $\mathscr{K}_s(\Iee_{0,\vartheta})$, and writes like in~(\ref{eq3-2-4}). Some basic considerations can be made on the order of this approximation. We know from~\citep{achieser} or~\citep{cheney} that, given $\varepsilon>0$, there exists $\eta$ such that for any $\mu_1,\mu_2\in[0,1]$, $|\mu_1-\mu_2|\leq \eta$ implies $|\Phi_\lambda(\mu_1)-\Phi_\lambda(\mu_2)|\leq\frac{\varepsilon}{2\vartheta}$, and that
$$
|\Phi_\lambda(\mu)-\Phi_{\lambda,\mathrm{app}}(\mu)|\leq\frac{\varepsilon}{2\vartheta}+\frac{\|\Phi_\lambda\|_{\mathscr{L}_\infty}}{2\vartheta\eta^2n}.
$$
There exists $\beta$ a positive bounded real number, such that $\left|\mu_1^{\frac{\lambda}{\alpha}}-\mu_2^{\frac{\lambda}{\alpha}}\right|\leq\beta|\mu_1-\mu_2|$, with for instance $\beta\geq\frac{\mathrm{e}^{\lambda\vartheta}-1}{1-\mathrm{e}^{-\alpha\vartheta}}$. This in turn implies that
$\eta\leq\frac{\varepsilon}{2\vartheta\beta}\mathrm{e}^{-2\lambda\vartheta}$. Taking the maximal admissible value for $\eta$, we finally obtain that the approximation order $n$ satisfies
$$
n\geq \frac{4\vartheta^3\beta^2}{\varepsilon^3}\mathrm{e}^{5\lambda\vartheta}. 
$$
From this simple consideration coming from the use of Bernstein polynomials, we obtain, to guarantee a norm upper bound $\varepsilon$ for the error, a condition on the order of the approximation. This condition is however quite conservative, and in practice, the order can be chosen iteratively, as illustrated in Fig.~\ref{fig4-1-1} and~Fig.~\ref{fig4-1-2}. 

%====================================================================
\subsection{Properties in time and frequency domains}\label{section4-2}
Let $f_{\mathrm{app}}$ be a kernel approximation in $\mathscr{B}(f,\varepsilon)$, for a given $\varepsilon>0$ and $f\in\mathscr{K}(\Iee_{\vartheta_1,\vartheta_2})$. The corresponding output is $y_{\mathrm{app}}(t)=(f_{\mathrm{app}}\ast u)(t)$. From~(\ref{eq2-1-5}), for any $1\leq p\leq\infty$ and $u\in\mathscr{L}_p(\Ree_+)$,
\begin{equation}\label{eq4-2-1}
\|y-y_{\mathrm{app}}\|_{\mathscr{L}_p}\leq \varepsilon\|u\|_{\mathscr{L}_p}.
\end{equation}
In other words, the output error $e(t)=y(t)-y_{\mathrm{app}}(t)$ can be made arbitrarily small for the $\mathscr{L}_p$-norm by a suitable choice of the arbitrary bound $\varepsilon$. Such a property includes the case of persistent inputs in $\mathscr{L}_1(\Ree_+)$. For the particular case $p=2$, we verify that the proposed approximation holds also in the frequency domain, since  
\begin{equation}\label{eq4-2-2}
\mathop{\mathrm{sup}}_{\mathrm{Re}\,s\geq 0}|\hat{f}(s)-\hat{f}_{\mathrm{app}}(s)|=\|\hat{f}-\hat{f}_{\mathrm{app}}\|_{\mathscr{H}_\infty}\leq\varepsilon.
\end{equation}

\begin{figure}[h!]
\begin{center}
\includegraphics[width=4.8in,height=3in]{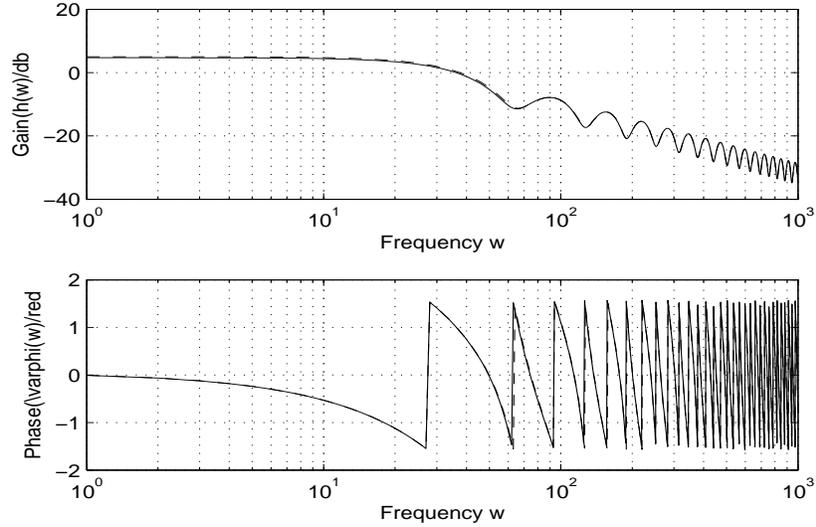}
\end{center}
\caption{Bode diagram of $\hat{\theta}_1(j\omega)$ and its kernel approximation $\hat{\theta}_{1,\mathrm{app}}(j\omega)$, with order $n=5$.}\label{fig4-2-2}
\end{figure}

\begin{thm}\label{def4-2-1}
Let $f\in\mathscr{A}$ be a given distribution, and $f_{\mathrm{app}}$ an approximation in $\mathscr{B}(f,\varepsilon)$, for some small $\varepsilon>0$. Then
$$
\|\hat{f}-\hat{f}_{\mathrm{app}}\|_{\mathscr{H}_\infty} =
\mathop{\mathrm{sup}}_{\omega\in\mathds{R}}{|\hat{f}(j\omega)-\hat{f}_{\mathrm{app}}(j\omega)|}
\leq\varepsilon,
$$
and for all $\omega$ in $\Ree$,
$$
|\mathrm{arg}(\hat{f}(j\omega))-\mathrm{arg}(\hat{f}_{\mathrm{app}}(j\omega))|\leq\varepsilon.
$$
\end{thm}
\begin{pf}
Since $f_{\mathrm{app}}\in\mathscr{B}(f,\varepsilon)$, and
$$
\|\hat{f}-\hat{f}_{\mathrm{app}}\|_{\mathscr{H}_\infty} 
\leq\|\hat{f}-\hat{f}_{\mathrm{app}}\|_{\hat{\mathscr{A}}},
$$
module inequality is trivial. To show that approximation holds also for the phase angle, we denote $\hat{f}_\mathrm{app}=\hat{f}+\hat{e}_{\mathrm{app}}$. Hence
$$
\left||\hat{f}(j\omega)|\,|1-\mathrm{e}^{j(\varphi(\omega)-\varphi_{\mathrm{app}}(\omega))}|-|\hat{e}_{\mathrm{app}}(j\omega)|\right| \leq \varepsilon,
$$
where $\varphi(\omega)=\mathrm{arg}(\hat{f}(j\omega))$ and $\varphi_{\mathrm{app}}(\omega)=\mathrm{arg}(\hat{f}_{\mathrm{app}}(j\omega))$.
This implies that
$$
|\hat{f}(j\omega)|\,|1-\mathrm{e}^{j(\varphi(\omega)-\varphi_{\mathrm{app}}(\omega))}|
$$
can be made as small as we please for all $\omega$, which in turn implies that $\varphi(\omega)$ and $\varphi_{\mathrm{app}}(\omega)$ are arbitrarily close for the $\mathscr{L}_\infty$-norm. 
\hfill{$\Box$}
\end{pf}
This frequency property is illustrated in Fig.~\ref{fig4-2-2}, where an approximation with order $n=5$ in $\mathscr{K}_s(\Iee_{0,1})$ of $\hat{\theta}_1(j\omega)$ is plotted.
The approximation of distributed delay in the graph topology over $\mathscr{A}$ yields to an approximation in both time and frequency domains. This is a strong property, that turns to be central in control problems.
%=====================================================================
\section{Application to control of time-delay systems}\label{section5}
Distributed delays appear naturally in the control of time-delay systems. We consider here as control application, the stabilization problem, and we illustrate the approximation method to realize the distributed time operator. 
 
\subsection{Stabilization}\label{section5-1}
\begin{figure}[h!]
\ifx\JPicScale\undefined\def\JPicScale{0.6}\fi
\psset{unit=\JPicScale mm}
\psset{linewidth=0.3,dotsep=1,hatchwidth=0.3,hatchsep=1.5,shadowsize=1,dimen=middle}
\psset{dotsize=0.7 2.5,dotscale=1 1,fillcolor=black}
\psset{arrowsize=1 2,arrowlength=1,arrowinset=0.25,tbarsize=0.7 5,bracketlength=0.15,rbracketlength=0.15}
\begin{pspicture}(-30,40)(160,70)
\pspolygon[](50,60)(70,60)(70,50)(50,50)
\psline{->}(5,55)(20,55)
\psline{->}(30,55)(50,55)
\psline{->}(70,55)(90,55)
\pspolygon[](120,60)(140,60)(140,50)(120,50)
\psline{->}(100,55)(120,55)
\psline{->}(95,70)(95,60)
\psline{->}(140,55)(160,55)
\psline(150,55)(150,42.5)
\psline(150,42.5)(25,42.5)
\rput(7.5,58.75){$u_1$}
\rput(145,20){}
\rput(60,60){}
\rput(60,55){}
\rput(130,55){}
\psline{->}(25,42.5)(25,50)
\rput{0}(25,55){\psellipse[](0,0)(5,-5)}
\rput{0}(95,55){\psellipse[](0,0)(5,-5)}
\rput(16.88,50){$+$}
\rput(20,45){$-$}
\rput(38.75,15){}
\rput(38.75,58.12){$e_{1}$}
\rput(108.75,58.75){$e_{2}$}
\rput(75.62,58.12){$y_{1}$}
\rput(91.25,70){$u_{2}$}
\rput(148.75,58.12){$y_{2}$}
\rput(100,65.62){$+$}
\rput(86.25,60){$+$}
\rput(59.38,55){$c$}
\rput(130,55){$p$}
\end{pspicture}
\caption{Feedback System.}\label{fig5-1-1}
\end{figure}
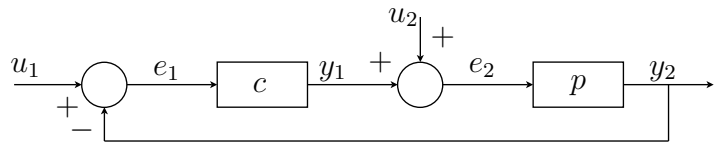

In the stabilization problem, we determine a control law of a given plant, such that the closed-loop system is stable. For this, consider
a distribution $p$ with a coprime factorization $(n,d)$ over $\mathscr{A}$, and define the closed-loop system described in Fig.~\ref{fig5-1-1}, with a compensator $c$ defined over the quotient field of $\mathscr{A}$. 

\noindent We assume here that the plant $p$ admits such a coprime factorization over $\mathscr{A}$. This is an assumption, but most dynamical systems fulfill it. The dynamical equations of this system are
\begin{equation}\label{eq5-1-1}
\left[\begin{matrix} y_1 \\ y_2 \end{matrix}\right] = 
H(p,c)\left[\begin{matrix} u_1 \\ u_2 \end{matrix}\right],\;\;
H(p,c)=\left[\begin{matrix} \frac{c}{1+pc} & \frac{-pc}{1+pc} \\ \frac{pc}{1+pc} & \frac{p}{1+pc} \end{matrix}\right].
\end{equation}
We say that $c$ stabilizes $p$, or the pair $(p,c)$ is stable, if 
the matrix $H(p,c)\in\mathscr{A}^{2\times 2}$. All internal signals in the closed-loop are bounded for any bounded exogeneous inputs $u_1$ and $u_2$. A necessary and sufficient condition for $(p,c)$ to be stable is the existence of some $n_c$ and $d_c$ in $\mathcal{A}$ such that
\begin{equation}\label{eq5-1-2}
n\ast n_c+d\ast d_c = \upsilon,
\end{equation}
with $\upsilon$ a unit in $\mathscr{A}$. Since $(n,d)$ is a coprime factorization,~(\ref{eq5-1-2}) holds, so we conclude that a stabilizing compensator $c$ has a coprime factorization $(n_c,d_c)$. A coprime factorization $(n_c,d_c)$ over $\mathscr{A}$ of a stabilizing compensator for $p$, includes, in general, some distributed delays. Approximating such a compensator as described in Section~\ref{section3} yields to the approximated controller $c_{\mathrm{app}}$ with factorization $(n_{c,\mathrm{app}},d_{c,\mathrm{app}})$, where $n_{c,\mathrm{app}}\in\mathscr{B}(n_c,\varepsilon_n)$ and $d_{c,\mathrm{app}}\in\mathscr{B}(d_c,\varepsilon_d)$, for some given $\varepsilon_n$ and $\varepsilon_d$. Applying this approximation in the control, we would know the conditions on $c_\mathrm{app}$, such that the pair $(p,c_{\mathrm{app}})$ is stable. The positive answer comes from a direct application of small gain theorem for BIBO-stability.

\begin{lem}\label{def5-1-1}
The pair $(p,c_{\mathrm{app}})$ is stable if 
$$
\mathrm{max}(\varepsilon_n,\varepsilon_d)<
\left\|\begin{pmatrix}n & d \end{pmatrix}
\right\|^{-1}_\mathscr{A}.
$$
\end{lem}
\begin{pf}
Let $(n_c,d_c)$ be a coprime factorization over $\mathscr{A}$ of a stabilizing controller. From~(\ref{eq5-1-2}), we have
$$
n\ast n_c+d\ast d_c = \delta.
$$
The approximated controller $c_{\mathrm{app}}$ with factorization $(n_{c,\mathrm{app}},d_{c,\mathrm{app}})$ yields in closed-loop
\begin{equation}\label{eq5-1-3}
n \ast n_{c,\mathrm{app}}+d\ast d_{c,\mathrm{app}} = \delta - \begin{pmatrix}n & d\end{pmatrix} \ast \begin{pmatrix}n_c-n_{c,\mathrm{app}} \\ d_c-d_{c,\mathrm{app}}\end{pmatrix}.
\end{equation}
The approximated controller will stabilize the plant $p$ if and only if the right hand side in~(\ref{eq5-1-3}) is a unit over $\mathscr{A}$. Since $\mathscr{A}$ is a Banach algebra, a sufficient condition is
$$
\left\| \begin{pmatrix}n & d\end{pmatrix} \ast \begin{pmatrix}n_c-n_{c,\mathrm{app}} \\ d_c-d_{c,\mathrm{app}}\end{pmatrix} \right\|_\mathscr{A}<1.
$$
Via approximations, $n_{c,\mathrm{app}}\in\mathscr{B}(n_c,\varepsilon_n)$ and $d_{c,\mathrm{app}}\in\mathscr{B}(d_c,\varepsilon_d)$, so that the above inequality yields the sufficient condition
$$
\left\| \begin{pmatrix}n & d\end{pmatrix} \right\|_\mathscr{A}\mathrm{max}(\varepsilon_n,\varepsilon_d)<1.
$$
\hfill{$\Box$}
\end{pf}
As a first comment, remark that it is always possible to determine an approximation $c_{\mathrm{app}}$ such that $(p,c_{\mathrm{app}})$ is stable. The counterpart will be in the order of approximation, that will increase when the required accuracy vanishes. Small gain theorem, which is still valid over any Banach algebra, gives us a sufficient condition on the approximation accuracy to guarantee robust stability. This condition helps us to determine the order of the approximation. Note that from~\citep{dahlehohta} where the conservativeness of small gain theorem for BIBO-stability was studied, a converse statement for Lemma~\ref{def5-1-1} holds, like in the $\mathscr{H}_\infty$ case. This highlights the weak conservation of such a condition for robust stabilization. As an application, consider the plant $\hat{y}(s) = \hat{p}(s) \hat{u}(s)$ given by
\begin{equation}\label{eq5-1-4}
\hat{p}(s) = \frac{\mathrm{e}^{-s}}{s-1}.
\end{equation}
A coprime factorization writes $\hat{n}=\frac{\mathrm{e}^{-s}}{s+1}$, $\hat{d}=\frac{s-1}{s+1}$, since
$$
\hat{n}(s)\,2\,\mathrm{e}^1+\hat{d}(s)\left(1+2\,\hat{\theta}_1(s)\right)=1.
$$
Hence, a stabilizing compensator for~(\ref{eq5-1-4}) is
$$
u(t) = -2(\theta_1\ast u)(t)+2\,\mathrm{e}^1\, y(t).
$$
From Lemma~\ref{def5-1-1}, a sufficient condition for robust stability with approximation of $\theta_1$ is that $\varepsilon_d\leq\frac{\mathrm{e}}{3+\mathrm{e}}$. From Fig.~\ref{fig4-1-2}, we see that a first order may be sufficient for stability purpose. In practice, a sufficient accuracy is obtained for a 5th order approximation over $\mathscr{K}_s(\Iee_{0,1})$, as can be seen in Fig.~\ref{fig4-1-1}, where step responses are plotted. Note also that if we modify the coprime factorization $(n,d)$, we could obtain lower order approximations.

\begin{figure}[h!]
\begin{center}
\includegraphics[width=4in,height=3in]{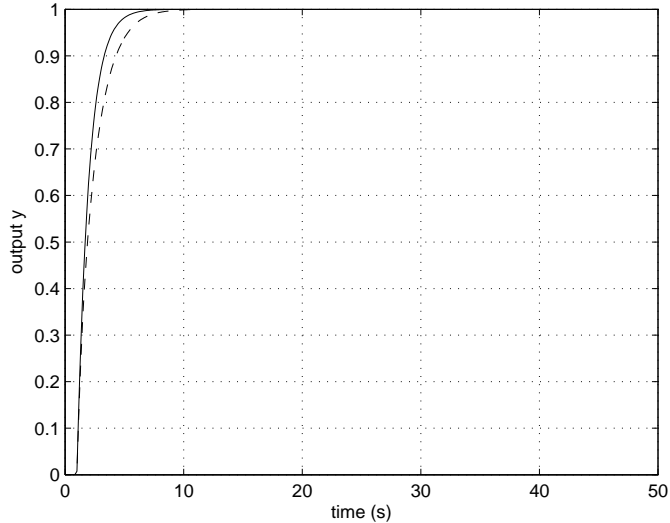}
\end{center}
\caption{Step responses for the desired closed-loop (continuous line) and for the closed-loop system with approximation of distributed delay (5th order, dash).}\label{fig5-2-2}
\end{figure}

%=====================================================================
\section{Conclusion}\label{section6}
In this paper, we proposed a general framework for rational approximation of distributed delay operators. The properties of this approximation were studied in both time and frequency domains. The effectiveness of this method was shown in simulation on the stabilization and the finite spectrum assignment problems, for general linear time-delay systems. We expect that this work will bring some new light in the understanding of distributed delay approximation, and more generally in approximation for the control of linear infinite dimensional systems. Such an approximation will provide foundations for a renewal of interest in control synthesis with distributed delays. Distributed delays inherit from integral control good robustness properties, and are at the core of numerous results in robustness analysis, optimization and control for time-delay systems. 


\bibliographystyle{model1-num-names}
\bibliography{HaoLuetal_SCL10}

\begin{thebibliography}{31}
\expandafter\ifx\csname natexlab\endcsname\relax\def\natexlab#1{#1}\fi
\providecommand{\bibinfo}[2]{#2}
\ifx\xfnm\relax \def\xfnm[#1]{\unskip,\space#1}\fi
%Type = Article
\bibitem[{Olbrot(1978)}]{olbrot}
\bibinfo{author}{A.~W. Olbrot},
\newblock \bibinfo{title}{Stabilizability, detectability, and spectrum
  assignment for linear autonomous systems with general time delays},
\newblock \bibinfo{journal}{IEEE Trans. on Autom. Contr.} \bibinfo{volume}{23}
  (\bibinfo{year}{1978}) \bibinfo{pages}{887--890}.
%Type = Article
\bibitem[{Kamen et~al.(1986)Kamen, Khargonekar, and Tannenbaum}]{kamenetal}
\bibinfo{author}{E.~W. Kamen}, \bibinfo{author}{P.~P. Khargonekar},
  \bibinfo{author}{A.~Tannenbaum},
\newblock \bibinfo{title}{Proper stable bezout factorizations and feedback
  control of linear time-delay systems},
\newblock \bibinfo{journal}{Int. J. Contr.} \bibinfo{volume}{43}
  (\bibinfo{year}{1986}) \bibinfo{pages}{837--857}.
%Type = Inproceedings
\bibitem[{Breth\'e and Loiseau(1996)}]{bretheloiseauFruit}
\bibinfo{author}{D.~Breth\'e}, \bibinfo{author}{J.~J. Loiseau},
\newblock \bibinfo{title}{A result that could bear fruit for the control of
  delay-differential systems},
\newblock in: \bibinfo{booktitle}{IEEE MSCA}, \bibinfo{address}{Chania,
  Greece}.
%Type = Article
\bibitem[{{Gl\"using-L\"uer{\ss}en}(1997)}]{glusingluerssen}
\bibinfo{author}{H.~{Gl\"using-L\"uer{\ss}en}},
\newblock \bibinfo{title}{A behavioral approach to delay-differential systems},
\newblock \bibinfo{journal}{SIAM J. Contr. Optimiz.} \bibinfo{volume}{35}
  (\bibinfo{year}{1997}) \bibinfo{pages}{480--499}.
%Type = Article
\bibitem[{Manitius and Olbrot(1979)}]{manitiusolbrot}
\bibinfo{author}{A.~Z. Manitius}, \bibinfo{author}{A.~W. Olbrot},
\newblock \bibinfo{title}{Finite spectrum assignment problem for systems with
  delays},
\newblock \bibinfo{journal}{IEEE Trans. on Autom. Contr.} \bibinfo{volume}{24}
  (\bibinfo{year}{1979}) \bibinfo{pages}{541--553}.
%Type = Article
\bibitem[{Smith(1959)}]{smith}
\bibinfo{author}{O.~J.~M. Smith},
\newblock \bibinfo{title}{A controller to overcome dead time},
\newblock \bibinfo{journal}{Inst. Soc. Amer. J.} \bibinfo{volume}{6}
  (\bibinfo{year}{1959}) \bibinfo{pages}{28--33}.
%Type = Article
\bibitem[{Artstein(1982)}]{artstein}
\bibinfo{author}{Z.~Artstein},
\newblock \bibinfo{title}{Linear systems with delayed controls: a reduction},
\newblock \bibinfo{journal}{IEEE Trans. on Autom. Contr.} \bibinfo{volume}{27}
  (\bibinfo{year}{1982}) \bibinfo{pages}{869--879}.
%Type = Article
\bibitem[{Dym et~al.(1995)Dym, Georgiou, and Smith}]{dymgeorgiousmithOptDelay}
\bibinfo{author}{H.~Dym}, \bibinfo{author}{T.~T. Georgiou},
  \bibinfo{author}{M.~C. Smith},
\newblock \bibinfo{title}{Explicit formulas for optimally robust controllers
  for delay systems},
\newblock \bibinfo{journal}{IEEE. Trans. on Autom. Contr.} \bibinfo{volume}{40}
  (\bibinfo{year}{1995}) \bibinfo{pages}{656--669}.
%Type = Inproceedings
\bibitem[{Mondi\'e et~al.(2001)Mondi\'e, Niculescu, and
  Loiseau}]{mondieniculesculoiseau}
\bibinfo{author}{S.~Mondi\'e}, \bibinfo{author}{S.~Niculescu},
  \bibinfo{author}{J.~J. Loiseau},
\newblock \bibinfo{title}{Delay robustness of closed loop finite assignment for
  input delay systems},
\newblock in: \bibinfo{booktitle}{IFAC Workshop on Time-Delay Systems}.
%Type = Inproceedings
\bibitem[{{Di Loreto}(2006)}]{diloreto}
\bibinfo{author}{M.~{Di Loreto}},
\newblock \bibinfo{title}{Finite time model matching for time-delay systems},
\newblock in: \bibinfo{booktitle}{IFAC Workshop on Time-Delay Systems},
  \bibinfo{address}{Aquila, Italy}.
%Type = Article
\bibitem[{Kamen et~al.(1985)Kamen, Khargonekar, and
  Tannenbaum}]{kamenetalfinite}
\bibinfo{author}{E.~W. Kamen}, \bibinfo{author}{P.~P. Khargonekar},
  \bibinfo{author}{A.~Tannenbaum},
\newblock \bibinfo{title}{Stabilization of time-delay systems using
  finite-dimensional compensators},
\newblock \bibinfo{journal}{IEEE Trans. on Autom. Contr.} \bibinfo{volume}{30}
  (\bibinfo{year}{1985}) \bibinfo{pages}{75--78}.
%Type = Article
\bibitem[{Partington(2004)}]{partington2004}
\bibinfo{author}{J.~R. Partington},
\newblock \bibinfo{title}{Some frequency-domain approaches to the model
  reduction of delay systems},
\newblock \bibinfo{journal}{Annual Reviews in Control} \bibinfo{volume}{28}
  (\bibinfo{year}{2004}) \bibinfo{pages}{65--73}.
%Type = Inproceedings
\bibitem[{{Van Assche} et~al.(1999){Van Assche}, Dambrine, and
  Lafay}]{vanassche}
\bibinfo{author}{V.~{Van Assche}}, \bibinfo{author}{M.~Dambrine},
  \bibinfo{author}{J.~Lafay},
\newblock \bibinfo{title}{Some problems arising in the implementation of
  distributed-delay control laws},
\newblock in: \bibinfo{booktitle}{38th IEEE Conference on Decision \& Control},
  pp. \bibinfo{pages}{4668--4672}.
%Type = Phdthesis
\bibitem[{Zhong(2003)}]{zhongphd}
\bibinfo{author}{Q.-C. Zhong}, \bibinfo{title}{Robust Control of Systems with
  Delays}, Ph.D. thesis, Imperial College London, \bibinfo{address}{London},
  \bibinfo{year}{2003}.
%Type = Inproceedings
\bibitem[{Santos and Mondi\'e(2000)}]{santosmondie}
\bibinfo{author}{O.~Santos}, \bibinfo{author}{S.~Mondi\'e},
\newblock \bibinfo{title}{Control laws involving distributed time delays:
  Robustness of the implementation},
\newblock in: \bibinfo{booktitle}{American Control Conference}, pp.
  \bibinfo{pages}{2479--2480}.
%Type = Article
\bibitem[{Mirkin(2004)}]{mirkinApproximation}
\bibinfo{author}{L.~Mirkin},
\newblock \bibinfo{title}{On the approximation of distributed-delay control
  laws},
\newblock \bibinfo{journal}{Systems \& Control Letters} \bibinfo{volume}{51}
  (\bibinfo{year}{2004}) \bibinfo{pages}{331--342}.
%Type = Article
\bibitem[{Mondi\'e and Michiels(2003)}]{mondiemichiels}
\bibinfo{author}{S.~Mondi\'e}, \bibinfo{author}{W.~Michiels},
\newblock \bibinfo{title}{Finite spectrum assignment of unstable time-delay
  systems with a safe implementation},
\newblock \bibinfo{journal}{IEEE Trans. on Autom. Contr.} \bibinfo{volume}{48}
  (\bibinfo{year}{2003}) \bibinfo{pages}{2207--2212}.
%Type = Inproceedings
\bibitem[{Zhong(2005)}]{zhong2005}
\bibinfo{author}{Q.-C. Zhong},
\newblock \bibinfo{title}{Rational implementation of distributed delay using
  extended bilinear transformations},
\newblock in: \bibinfo{booktitle}{IFAC World Congress},
  \bibinfo{address}{Prague, Czech Republic}.
%Type = Article
\bibitem[{Richard(2003)}]{richard}
\bibinfo{author}{J.-P. Richard},
\newblock \bibinfo{title}{Time-delay systems: An overview of some recent
  advances and open problems},
\newblock \bibinfo{journal}{Automatica} \bibinfo{volume}{39}
  (\bibinfo{year}{2003}) \bibinfo{pages}{1667--1694}.
%Type = Article
\bibitem[{Zhong(2004)}]{zhong2004}
\bibinfo{author}{Q.-C. Zhong},
\newblock \bibinfo{title}{On distributed delay in linear control laws. part 1:
  Discrete-delay implementation},
\newblock \bibinfo{journal}{IEEE Trans. on Autom. Contr.} \bibinfo{volume}{49}
  (\bibinfo{year}{2004}) \bibinfo{pages}{2074--2080}.
%Type = Book
\bibitem[{Achieser(1956)}]{achieser}
\bibinfo{author}{N.~I. Achieser}, \bibinfo{title}{Theory of approximation},
  \bibinfo{publisher}{Frederick Ungar Publishing Corp., New York},
  \bibinfo{year}{1956}.
%Type = Book
\bibitem[{Cheney(1982)}]{cheney}
\bibinfo{author}{E.~W. Cheney}, \bibinfo{title}{Introduction to approximation
  theory}, \bibinfo{publisher}{AMS, 2nd Ed., Chelsea}, \bibinfo{year}{1982}.
%Type = Article
\bibitem[{Kammler(1976)}]{kammler}
\bibinfo{author}{D.~W. Kammler},
\newblock \bibinfo{title}{Approximation with sums of exponentials in
  $l_p[0,\infty)$},
\newblock \bibinfo{journal}{J. of Approximation Theory} \bibinfo{volume}{16}
  (\bibinfo{year}{1976}) \bibinfo{pages}{384--408}.
%Type = Article
\bibitem[{Vidyasagar and Anderson(1989)}]{vidyasagaranderson}
\bibinfo{author}{M.~Vidyasagar}, \bibinfo{author}{B.~D.~O. Anderson},
\newblock \bibinfo{title}{Approximation and stabilization of distributed
  systems by lumped systems},
\newblock \bibinfo{journal}{Systems \& Control Letters} \bibinfo{volume}{12}
  (\bibinfo{year}{1989}) \bibinfo{pages}{95--101}.
%Type = Article
\bibitem[{Otha et~al.(1992)Otha, Maeda, and Kodama}]{ohtaetal}
\bibinfo{author}{Y.~Otha}, \bibinfo{author}{H.~Maeda},
  \bibinfo{author}{S.~Kodama},
\newblock \bibinfo{title}{Rational approximation of $l_{1}$ optimal controllers
  for siso systems},
\newblock \bibinfo{journal}{IEEE Trans. on Autom. Contr.} \bibinfo{volume}{37}
  (\bibinfo{year}{1992}) \bibinfo{pages}{1683--1691}.
%Type = Article
\bibitem[{Callier and Desoer(1978)}]{callierdesoer}
\bibinfo{author}{F.~M. Callier}, \bibinfo{author}{C.~A. Desoer},
\newblock \bibinfo{title}{An algebra of transfer functions for distributed
  linear time-invariant systems},
\newblock \bibinfo{journal}{IEEE Trans. Circuits Syst.} \bibinfo{volume}{25}
  (\bibinfo{year}{1978}) \bibinfo{pages}{651--662}.
%Type = Book
\bibitem[{Desoer and Vidyasagar(1975)}]{desoervidyasagar}
\bibinfo{author}{C.~A. Desoer}, \bibinfo{author}{M.~Vidyasagar},
  \bibinfo{title}{Feedback systems: Input-output properties},
  \bibinfo{publisher}{Academic Press}, \bibinfo{address}{New-York},
  \bibinfo{year}{1975}.
%Type = Article
\bibitem[{Breth\'e and Loiseau(1997)}]{bretheloiseauFSA}
\bibinfo{author}{D.~Breth\'e}, \bibinfo{author}{J.~J. Loiseau},
\newblock \bibinfo{title}{An effective algorithm for finite spectrum assignment
  of single-input systems with delays},
\newblock \bibinfo{journal}{J. Math. Computer in Simulations}
  \bibinfo{volume}{Special issue on time-delay systems, V. Kolmanovskii and J.
  P. Richard Eds.} (\bibinfo{year}{1997}).
%Type = Book
\bibitem[{Vidyasagar(1985)}]{vidyasagarbook}
\bibinfo{author}{M.~Vidyasagar}, \bibinfo{title}{Control System Synthesis. A
  Factorization Approach}, \bibinfo{publisher}{MIT Press},
  \bibinfo{address}{Cambridge, Massachussets}, \bibinfo{year}{1985}.
%Type = Book
\bibitem[{Ralston and Rabinowitz(2001)}]{ralstonrabinowitz}
\bibinfo{author}{A.~Ralston}, \bibinfo{author}{P.~Rabinowitz},
  \bibinfo{title}{A first course in numerical analysis},
  \bibinfo{publisher}{Dover Publications, New York, 2nd Ed.},
  \bibinfo{year}{2001}.
%Type = Article
\bibitem[{Dahleh and Otha(1988)}]{dahlehohta}
\bibinfo{author}{M.~A. Dahleh}, \bibinfo{author}{Y.~Otha},
\newblock \bibinfo{title}{A necessary and sufficient condition for robust bibo
  stability},
\newblock \bibinfo{journal}{Systems \& Control Letters} \bibinfo{volume}{11}
  (\bibinfo{year}{1988}) \bibinfo{pages}{271--275}.

\end{thebibliography}

\end{document}